\documentclass[final,5p,twocolumn,preprint,nopreprintline,times,authoryear]{elsarticle}

\usepackage{hyperref}
\usepackage{amsmath}
\usepackage{amssymb}
\usepackage{cuted}
\usepackage{graphicx}
\usepackage{caption}
\usepackage{subcaption}

\graphicspath{
    {plots/}
}

\journal{Journal of Mathematical Analysis and Applications}

\begin{document}

\begin{frontmatter}

\title{Approximating the mode of the non-central chi-squared distribution}

\author[uio]{Victor Ananyev}
\author[uio]{Alexander Lincoln Read}
\affiliation[uio]{organization={Department of Physics, University of Oslo},
            addressline={Blindern},
            city={Oslo},
            postcode={0316},
            state={Oslo},
            country={Norway}}

\begin{abstract}
In this paper we consider the probability density function (PDF) of the non-central $\chi^2$ distribution with arbitrary number of degrees of freedom and non-centrality. For this function we find the approximate location of the maximum and discuss related edge cases of 1 and 2 degrees of freedom. We also use this expression to demonstrate the improved performance of the C++ Boost's implementation of the non-central $\chi^2$ and extend the domain of its applicability.
\end{abstract}



\end{frontmatter}


\section{Introduction}\label{introduction}
Properties of the non-central $\chi^2$ distribution were described before in literature~(\cite{Andrs2008, hogan2013, Saulis2001}). However, the topic of the mode of the non-central $\chi^2$ was significantly underrepresented. We would like to focus on the mode specifically in this paper.

Let $X_1, X_2, ..., X_n$ be normally distributed random variables with unit variance and means $\mu_1, \mu_2, ..., \mu_n$. The sum $X_1^2 + X_2^2 + ... + X_n^2$ follows the non-central $\chi^2$ distribution with $k = n$ degrees of freedom and non-centrality $\lambda = \mu_1^2 + \mu_2^2 + ... + \mu_n^2$. The probability density function of this distribution has a closed form expression:
\begin{align}
    &f_{k, \lambda}(x) = \frac{1}{2} \exp^{-\frac{x + \lambda}{2}} \left(\frac{x}{\lambda}\right)^{\frac{k-2}{4}} I_{\frac{k-2}{2}}(\sqrt{\lambda x})\;,
\end{align}
where $I_{\nu}(x)$ is a modified Bessel function of the first kind.

\begin{figure}[h]
	\centering
	\begin{subfigure}[t]{0.22\textwidth}
		\centering
		\includegraphics[width=\linewidth]{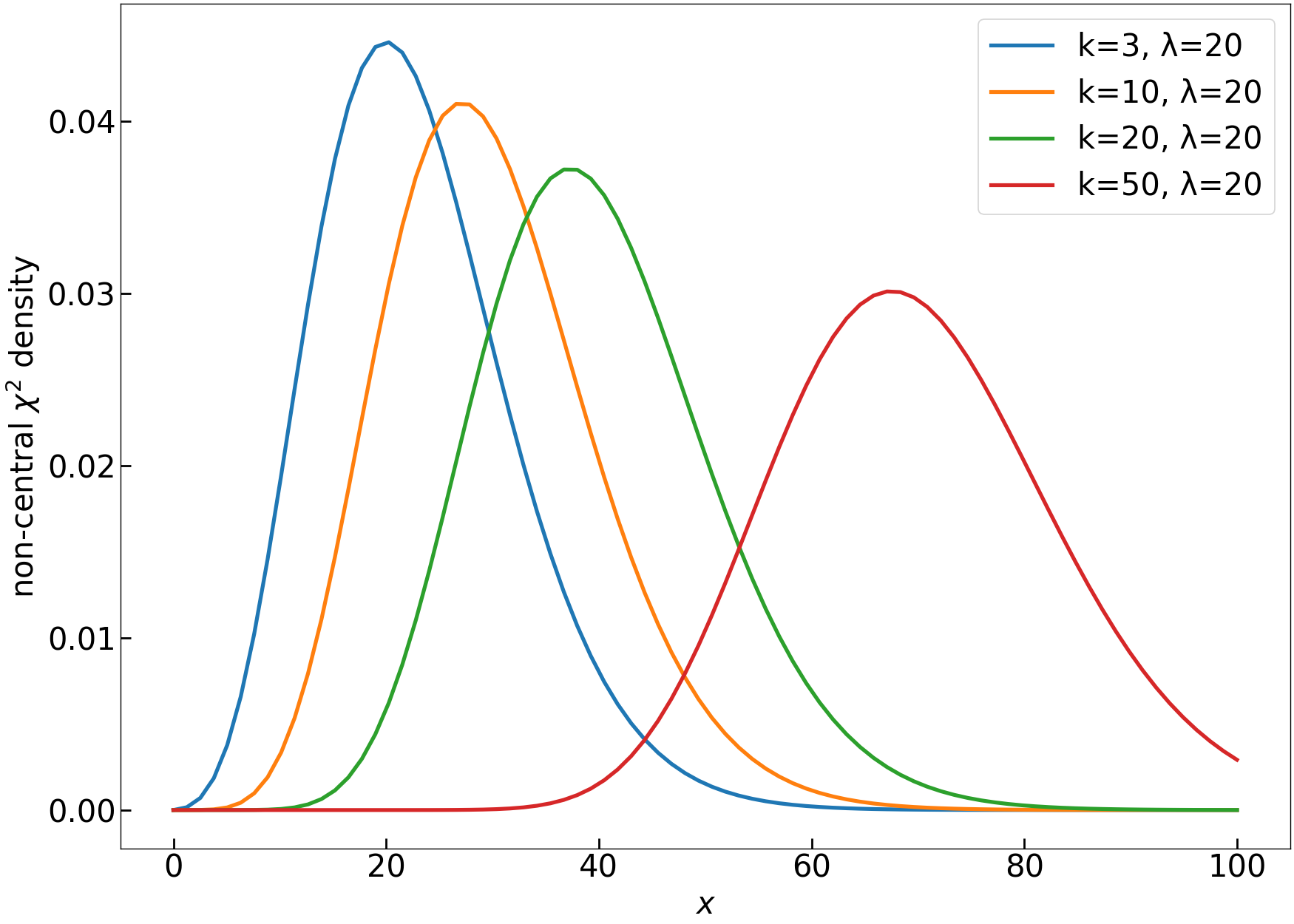}
		\caption{Dependency on the number of d.o.f.~$k$.}\label{fig:ncx2-curves-df}
	\end{subfigure}
	\quad
	\begin{subfigure}[t]{0.22\textwidth}
		\centering
		\includegraphics[width=\linewidth]{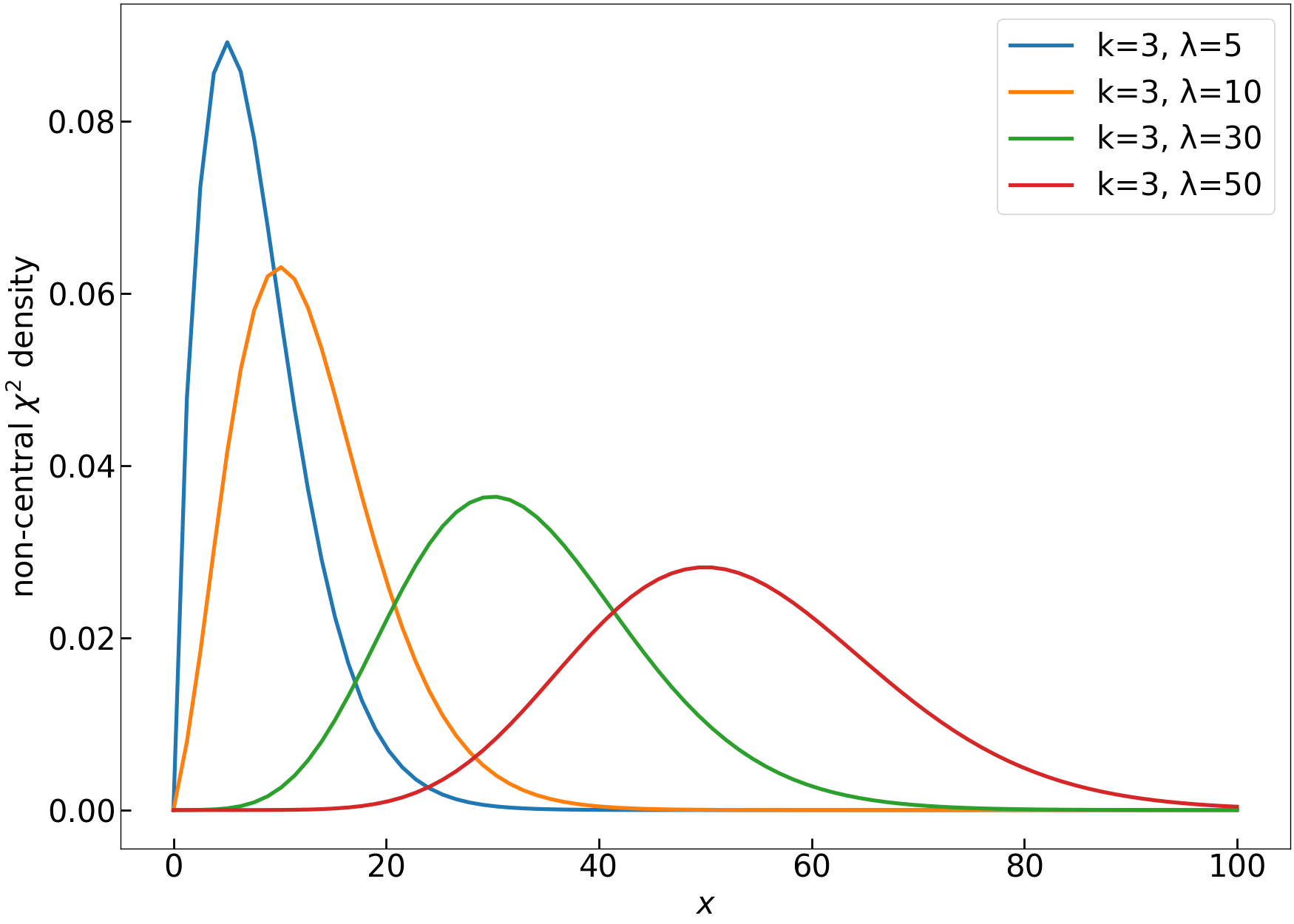}
		\caption{Dependency on the non-centrality $\lambda$.}\label{fig:ncx2-curves-lambda}
	\end{subfigure}
	\caption{Non-central $\chi^2$ distributions and behavior of the mode.}\label{fig:ncx2-curves}
\end{figure}

We are interested in the value of $x_{mode}$ that maximizes $f_{k, \lambda}(x)$. Typical shapes of the pdf of the non-central $\chi^2$ distribution are shown in Fig.~\ref{fig:ncx2-curves}.

\begin{figure}[h]
	\centering
 	\includegraphics[width=0.44\textwidth]{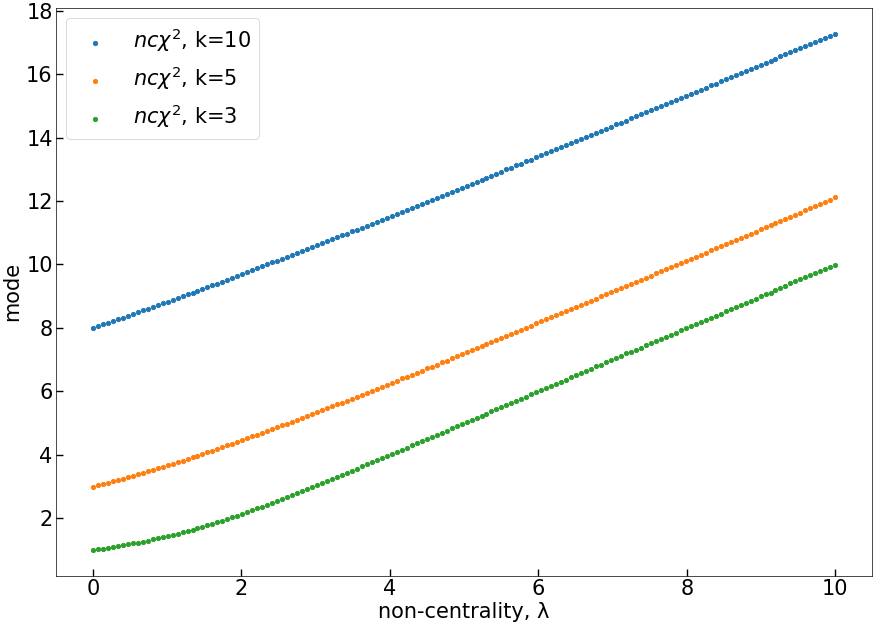}
	\caption{Mode of the non-central $\chi^2$ as a function of the non-centrality parameter $\lambda$.}\label{fig:ncx2-max}
\end{figure}

When the number of degrees of freedom $k$ is fixed, we can plot the dependency of the maximum of the pdf as a function of the non-centrality parameter $\lambda$, see Fig.~\ref{fig:ncx2-max}.

We observe that the bigger $\lambda$ is the better the mode appears to be approximated with a straight line. The derivation of the line parameters together with the analysis of the edge cases of small number of degrees of freedom, where the mode does not exist, constitute the main results of the paper.

\section{Derivation}

\subsection{Master equation}
In this section we obtain the transcendental equation (Eq.~\ref{eq:master}) that determines the mode of the non-central $\chi^2$ distribution. We reduce it to the ordinary differential equation (Eq.~\ref{eq:master-ode}), where the non-centrality parameter $\lambda$ is the argument, and the number of degrees of freedom $k$ is a parameter. Finally, we solve the ODE approximately with a Taylor expansion (Eq.~\ref{eq:taylor-coef-sols}) and investigate edge cases of 1 and 2 degrees of freedom (Sec.~\ref{subsec:edge-cases}).

We start by setting the derivative of the density of the non-central $\chi^2$ (Eq.~\ref{eq:ncx2-pdf-to-zero}) to zero. This leads us to the transcendental equation (Eq.~\ref{eq:master}) that determines the mode of the distribution:
\begin{align}
    &\frac{d}{dx} \chi^2_{k, \lambda}(x) = \frac{1}{2} \chi^2_{d, \lambda}(x) \cdot \left[ -1 + \frac{k-2}{2 x} + \sqrt{\frac{\lambda}{x}}\frac{I^{\prime}_{\frac{k-2}{2}}(\sqrt{\lambda x})}{I_{\frac{k-2}{2}}(\sqrt{\lambda x})} \right]\label{eq:ncx2-pdf-to-zero}\;, \\
    &\frac{d}{dx} \chi^2_{k, \lambda}(x) = 0 \Rightarrow \sqrt{\lambda x} I^{\prime}_{\frac{k-2}{2}}(\sqrt{\lambda x}) = (x - \frac{k-2}{2}) I_{\frac{k-2}{2}}(\sqrt{\lambda x})\label{eq:master}\;.
\end{align}
We can eliminate the derivative in Eq.~\ref{eq:master} by using the differential equation for the modified Bessel function~\cite[Eq.~10.25.1]{NIST:DLMF}:
\begin{align}
    &t^2 \frac{d^2}{dt^2} I_{\nu}(t) + t \frac{d}{dt} I_{\nu}(t) - (t^2 + \nu^2) I_{d, \lambda}(t) = 0 \label{eq:bessel-definition}\;.
\end{align}
To make use of Eq.~\ref{eq:bessel-definition}, we need the expression for $I^{\prime \prime}_{\frac{k-2}{2}}$, therefore, we differentiate Eq.~\ref{eq:master} by $\lambda$. Since the mode depends on the non-centrality $\lambda$, we should remember that $x = x(\lambda)$, thus $\frac{dx}{d\lambda} = x^{\prime}$. The resulting expression for $I^{\prime \prime}_{\frac{k-2}{2}}$ is as follows:
\begin{align}
    &\sqrt{\lambda x} I^{\prime \prime}_{\frac{k-2}{2}} (\sqrt{\lambda x}) =
    (x - \frac{k}{2}) I^{\prime}_{\frac{k-2}{2}}(\sqrt{\lambda x}) +  \frac{2 \sqrt{\lambda x} x^{\prime}}{x + \lambda x^\prime} I_{\frac{k-2}{2}}(\sqrt{\lambda x})\label{eq:master-derivative}\;.
\end{align}
We substitute $I^{\prime}_{\frac{k-2}{2}}$ (Eq.~\ref{eq:master}) and $I^{\prime \prime}_{\frac{k-2}{2}}$ (Eq.~\ref{eq:master-derivative}) into the differential equation for the modified Bessel function (Eq.~\ref{eq:bessel-definition}). We then use~\cite[Eq.~10.29.4]{NIST:DLMF} to decrease the order of the derivatives of the modified Bessel functions. Assuming that the Bessel function itself is non-zero at the mode, we arrive to the following differential equation for the mode as a function of the non-centrality parameter $\lambda$:
\begin{align}
    &\lambda x^{\prime} (x - k - \lambda + 4) + x (x - k - \lambda + 2) = 0\;.\label{eq:master-ode}
\end{align}
\subsection{Approximate solution}
\begin{figure}[h]
	\centering
 	\includegraphics[width=0.44\textwidth]{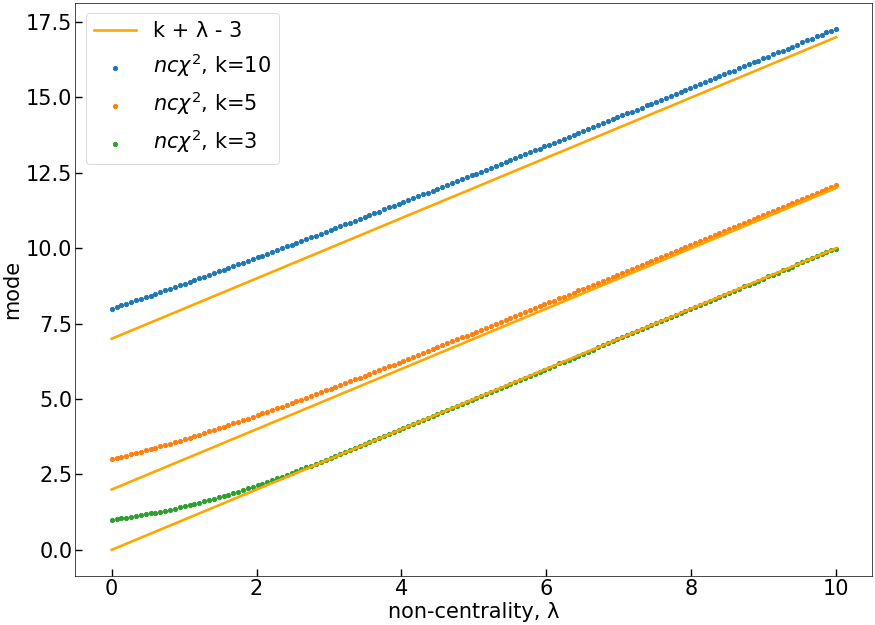}
	\caption{Linear approximation to the mode of the non-central $\chi^2$ compared to the more precise numerical solution as a function of the non-centrality parameter $\lambda$.}\label{fig:ncx2-vs-approx}
\end{figure}

We observed that the linear approximation works better with growing $\lambda$, thus we introduce the asymptotic parameter $t = \frac{k}{\lambda} << 1$ to build the expansion. We expect the solution to be linear in $\lambda$, however the asymptotic expansion of $x(t) = C_0 + C_1 t + ...$ won't provide us with a solution linear in $\lambda$. Therefore, we reparametrize $x(t)$ with a new function $y(t) = t x(t)$:
\begin{align}
    &t = \frac{k}{\lambda}\;,\label{eq:asymp-scale} \\
    &y(t) = t x(t)\;.
\end{align}
We obtain the following equation after the reparametrization:
\begin{align}
   &-(y' t - y) (y - kt - k + 4t) + y (y - kt - k + 2t) = 0\label{eq:reduced-ready-eq}\;.
\end{align}
To solve Eq.~\ref{eq:reduced-ready-eq}, we expand $y(t)$ into the Taylor series by the scale parameter $t = \frac{k}{\lambda}$. We would like to find the linear solution and one extra term that estimates the error. Thus, we cut the series at the third power of $t$ in order to account for the derivative. After solving algebraic equations for the coefficients near each power of $t$, we arrive to the resulting approximate expression for the mode:
\begin{align}
    &y(t) = C_0 + C_1 t + C_2 t^2 + C_3 t^3 + O(t^4)\;, \\
    &C_0 = k, \quad C_1 = k-3, \quad C_2 = \frac{k-3}{2k}\;,\label{eq:taylor-coef-sols}\\
    &\boxed{x_{mode} = \lambda + k - 3 + \frac{k-3}{2 \lambda} + O\left(\frac{k^2}{\lambda^2}\right)} \label{eq:main-result}\;.
\end{align}
We plot the linear approximation Eq.~\ref{eq:main-result} together with the precise numerical solution Fig.~\ref{fig:ncx2-max} in order to verify the approximation is correct, see Fig.~\ref{fig:ncx2-vs-approx}.

\subsection{Small number of degrees of freedom}\label{subsec:edge-cases}
\paragraph{Case $k < 2$}
The asymptotic behavior of the modified Bessel function at $x \rightarrow 0$~\cite[Eq.~10.30.1]{NIST:DLMF} shows that the PDF of the non-central $\chi^2$ diverges, thus it doesn't have a mode:
\begin{align}
 &\chi^2_{k, \lambda}(x) \rightarrow \frac{1}{2\Gamma(\frac{k}{2})} \frac{1}{(2 \lambda)^{\frac{k - 2}{2}}} \mathrm{e}^{-\frac{\lambda}{2}} \left(\sqrt{\lambda x}\right)^{k-2}, \quad x \rightarrow 0\;.
\end{align}
\paragraph{Case $k = 2$}
In this case, the PDF at $x=0$ is finite. If the derivative at $x=0$ is positive, then the maximum is not there. The expression for the derivative (Eq.~\ref{eq:ncx2-pdf-k2}) and its asymptotic behavior at $x \rightarrow 0$ (Eq.~\ref{eq:ncx2-pdf-k2-0}) are shown below:
\begin{align}
 &\frac{d}{dx} \chi^2_{k, \lambda}(x) = \frac{1}{2} \chi^2_{d, \lambda}(x) \cdot \left[ -1 +  \sqrt{\frac{\lambda}{x}}\frac{I_{-1}(\sqrt{\lambda x})}{I_{0}(\sqrt{\lambda x})} \right]\label{eq:ncx2-pdf-k2}\;, \\
 &\frac{d}{dx} \chi^2_{k, \lambda}(x) \rightarrow \frac{1}{2} \chi^2_{d, \lambda}(x) \cdot \left[ -1 + \frac{\lambda}{2} \right], \quad x \rightarrow 0\label{eq:ncx2-pdf-k2-0}\;.
\end{align}
We observe that when $\lambda > 2$, the PDF of the non-central $\chi^2$ doesn't have its maximum at $x = 0$. In the region $\lambda < 2$, the asymptotic scale $t = \frac{k}{\lambda} > 1$, hence our approximation is inapplicable in this region and we refrain from analysing it.

\section{Application}
There exist a number of numerical procedures for finding the mode of a distribution~\cite[Ch.~10]{NR}. Some of them require the search region to be specified.

For example, the widely used C++ library Boost~(\cite{boost}) identifies the search region based on an initial guess for the mode $x_0$. Boost iteratively checks regions of the form $[x_0/2, 2x_0]$,  $[x_0/2^2, 2^2x_0]$, etc. When the value of the PDF at both ends of the region becomes smaller than the value at the initial guess point $x_0$, the algorithm initiates the search for the maximum inside of the region.

At the time of writing, Boost used $x_0 = k + 1$ as the initial guess. We already know, based on the approximate solution (Eq.~\ref{eq:main-result}), that the chosen guess will undershoot at large non-centrality values $\lambda$. Let's estimate $\lambda$ above which the method will require the second iteration for the region to cover the mode. For this we compare the linear estimate for the location of the mode~(Eq.~\ref{eq:main-result}) to the initial guess $x_0$ used by Boost:
\begin{align}
    &k + \lambda - 3 > 2 \cdot (k + 1)\;, \\
    &\lambda > k + 5 \label{eq:boost-lam-thresh-simple}\;.
\end{align}
With Eq.~\ref{eq:boost-lam-thresh-simple}, for any number of d.o.f. $k$ we are able to specify the threshold $\alpha$, defined by $\lambda  = \alpha k$, at which the original initial guess starts undershooting:
\begin{align}
    &\alpha > \frac{5}{k} + 1 \label{eq:boost-alpha-thresh}\;.
\end{align}
We see that large $k$ corresponds to small thresholds $\alpha$. The most conservative estimate for the threshold would be at the smallest $k$ possible: $k = 2$. Thus, $\alpha = 3.5$ is the threshold that approximately works for $k=2$ and is the overestimated threshold for the bigger values of $k$.

The threshold $\alpha$~(Eq.~\ref{eq:boost-alpha-thresh}) is closely related to the asymptotic scale $t$~(Eq.~\ref{eq:asymp-scale}) that we used for finding the approximate solution, specifically: $t = \frac{k}{\lambda} = \frac{1}{\alpha}$. For example, the conservative threshold $\alpha \approx 3.5$ corresponds to the asymptotic scale $t \approx 0.25 < 1$. It means that the region where the original guess of Boost undershoots, is, at the same time, the region where our approximate solution for the mode becomes applicable and can be used as a corrected initial guess. However, the fact that we use the conservative threshold may lead to the situation where the original method has already started undershooting but $\lambda$ is not yet big enough to turn on the corrected regime.

\begin{figure}[h]
	\centering
	\begin{subfigure}[t]{0.35\textwidth}
		\includegraphics[width=0.92\linewidth]{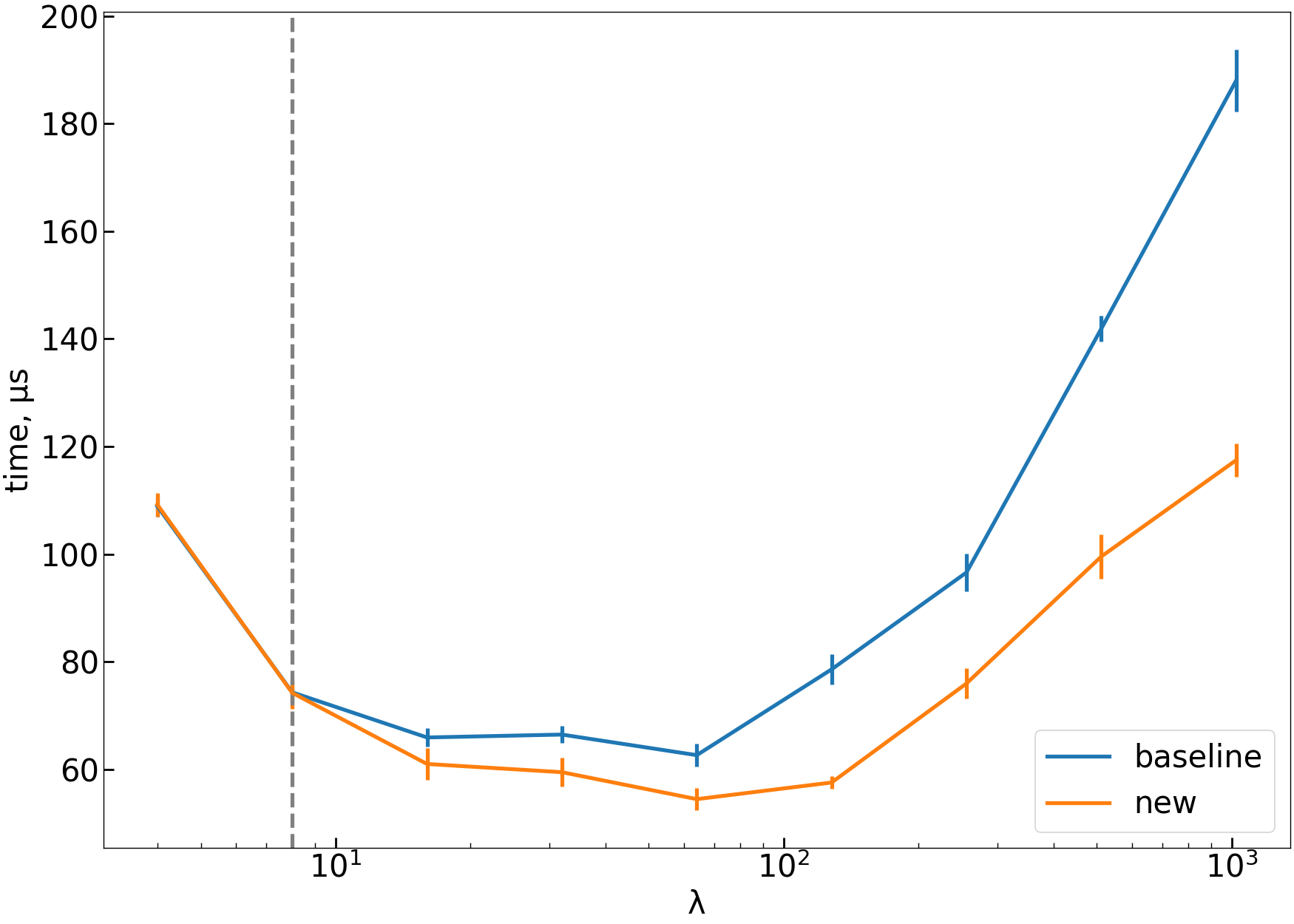}
        \caption{}\label{boost-k-fixed_2}
	\end{subfigure}
	\;
	\begin{subfigure}[t]{0.35\textwidth}
		\includegraphics[width=0.92\linewidth]{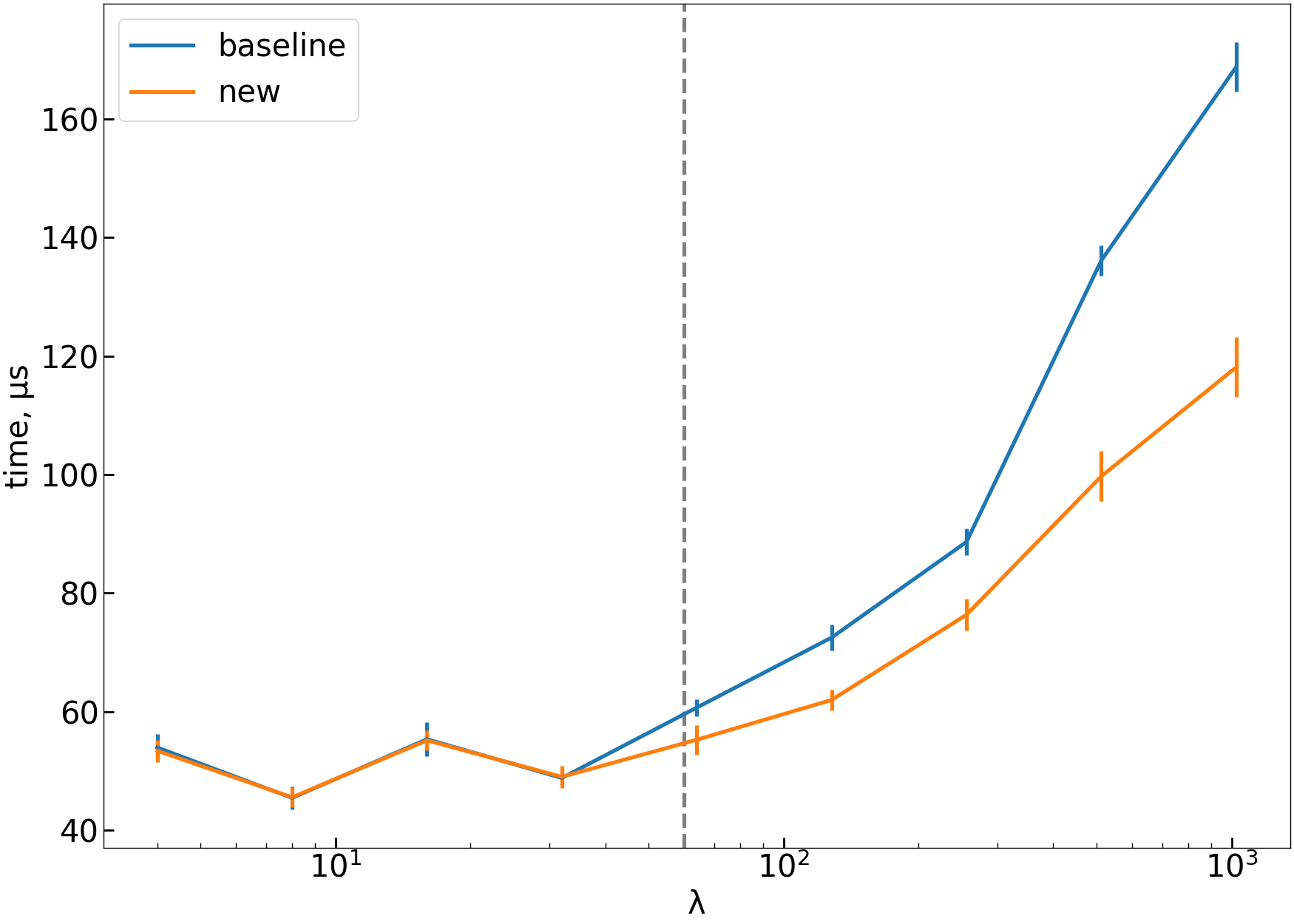}
		\caption{}\label{boost-k-fixed_15}
	\end{subfigure}
	\;
	\begin{subfigure}[t]{0.35\textwidth}
		\includegraphics[width=0.92\linewidth]{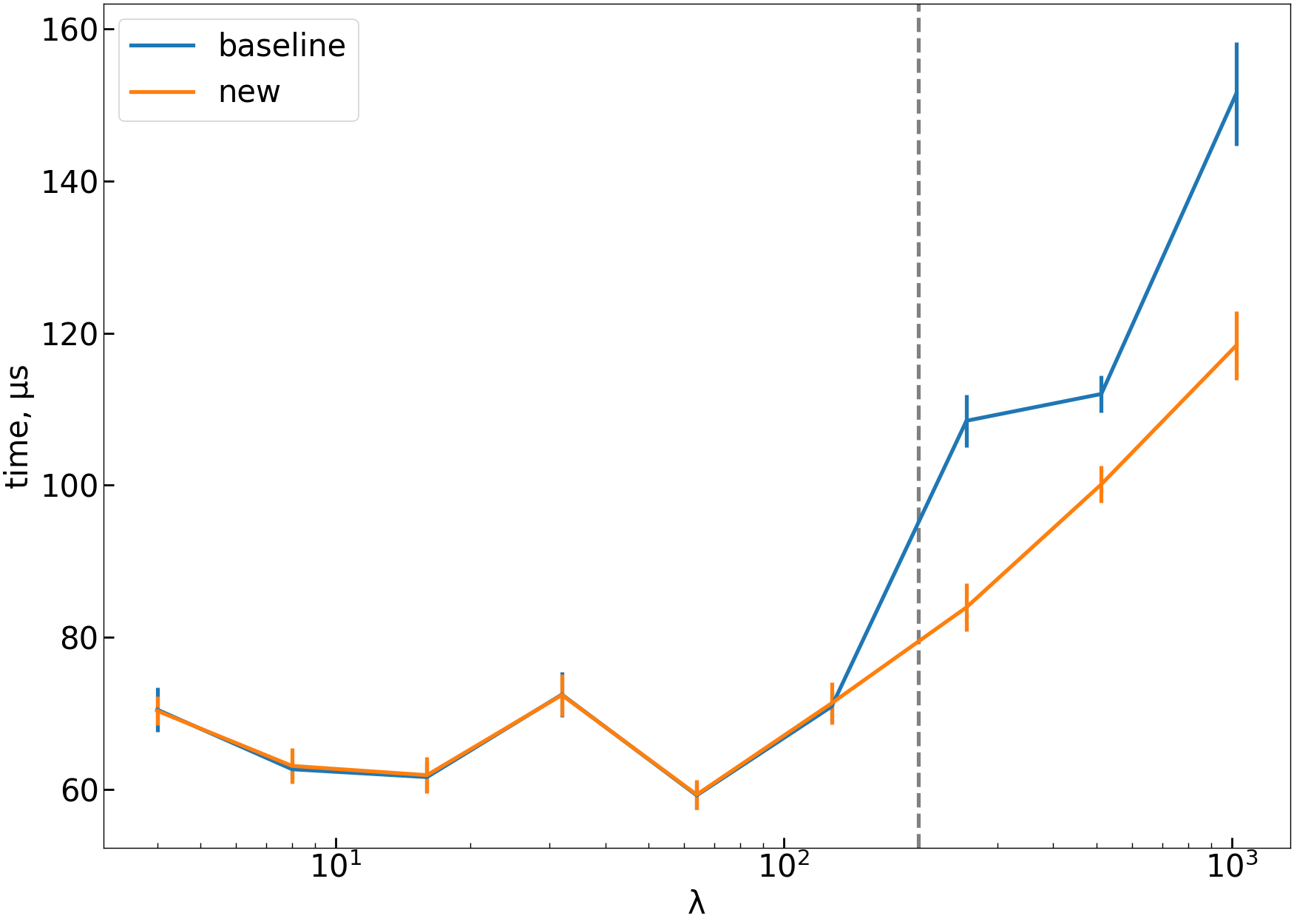}
		\caption{}\label{boost-k-fixed_50}
	\end{subfigure}
    \caption{Run time as a function of the non-centrality $\lambda$ for d.o.f. $k=2$~(\ref{boost-k-fixed_2}), $k=15$~(\ref{boost-k-fixed_15}), $k=50$~(\ref{boost-k-fixed_50}). Vertical line shows the threshold at which the corrected expression replaces the original initial guess.}\label{fig:boost-k-fixed}
\end{figure}

\begin{figure}[t]
    \centering
	\begin{subfigure}[t]{0.46\textwidth}
        \includegraphics[width=0.93\linewidth]{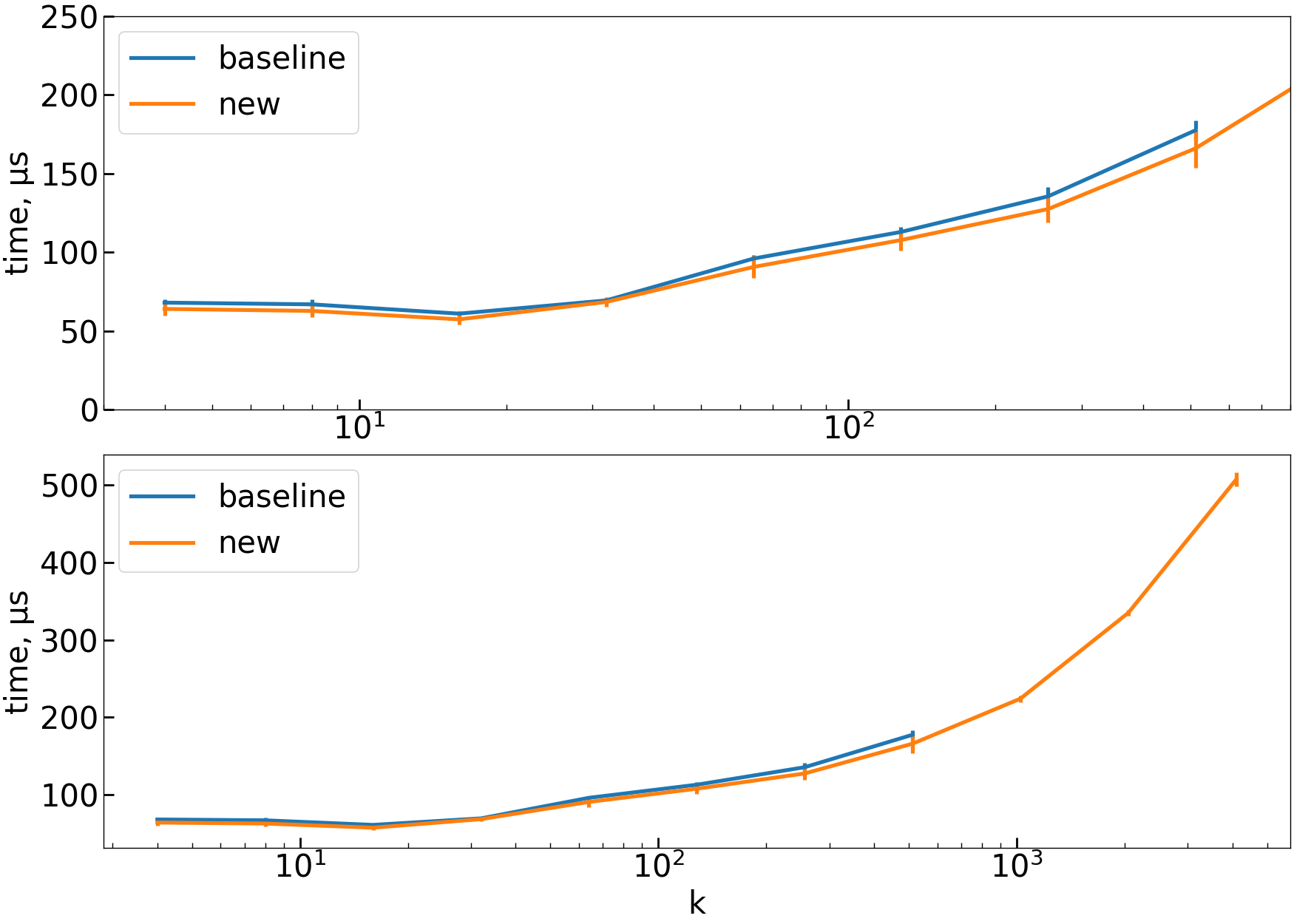}
        \caption{}\label{boost-scale-fixed_0_25}
	\end{subfigure}
	\quad
	\begin{subfigure}[t]{0.46\textwidth}
		\includegraphics[width=0.93\linewidth]{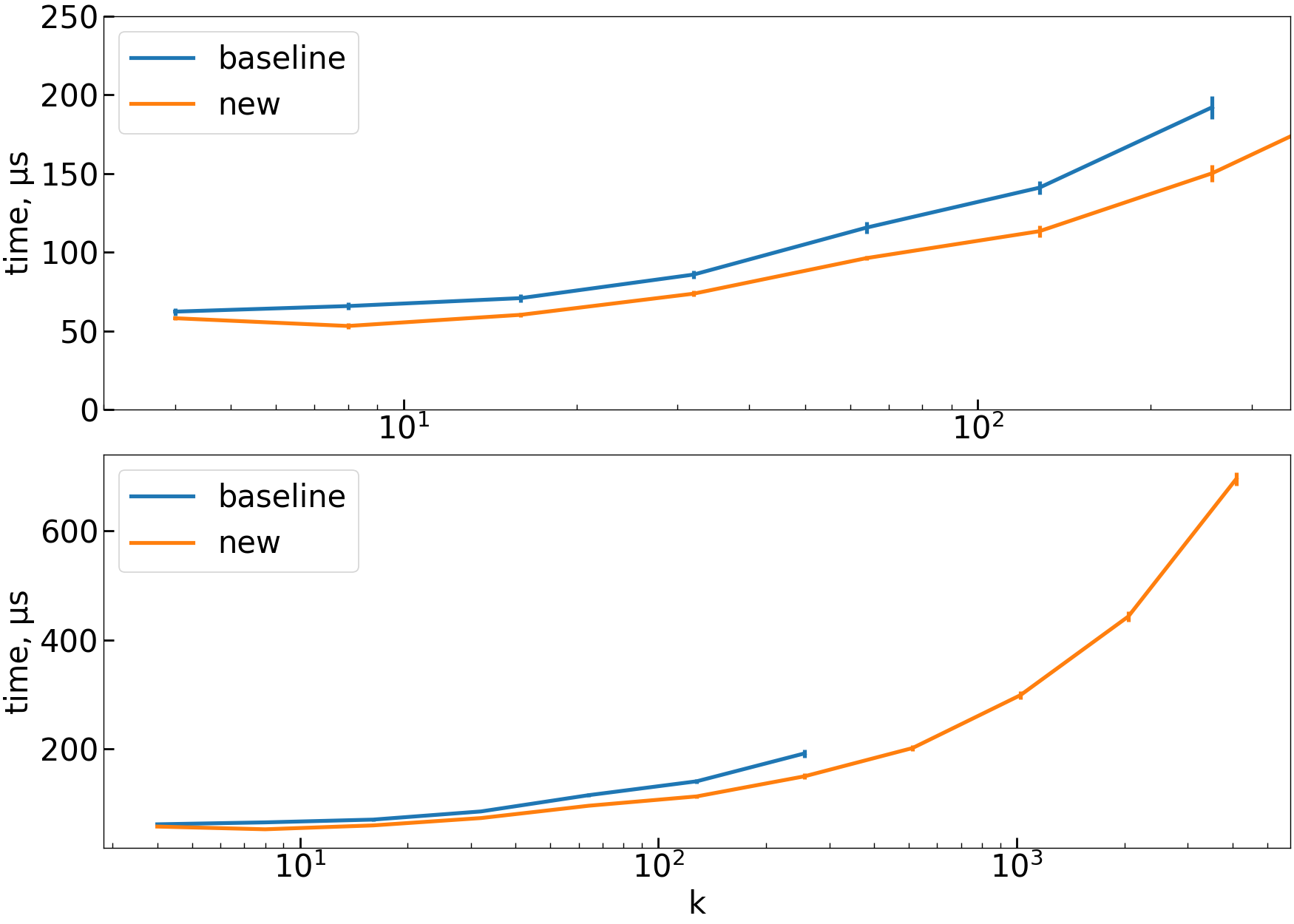}
		\caption{}\label{boost-scale-fixed_0_15}
	\end{subfigure} \\
	\begin{subfigure}[t]{0.46\textwidth}
		\includegraphics[width=0.93\linewidth]{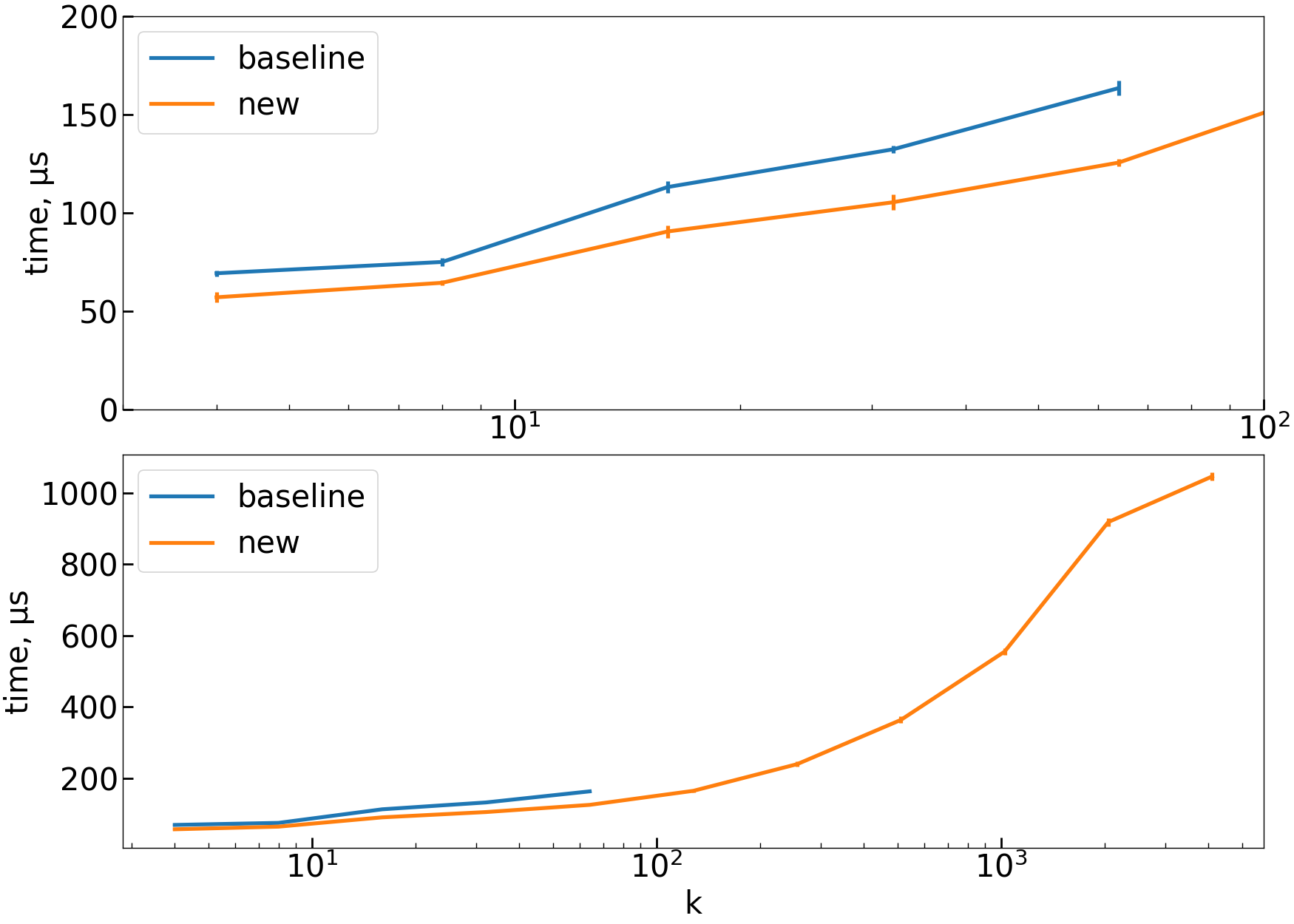}
		\caption{}\label{boost-scale-fixed_0_05}
	\end{subfigure}
    \caption{Run time as a function of the number of d.o.f. ($k$) for the asymptotic scale values $\frac{k}{\lambda}=0.25$~(\ref{boost-scale-fixed_0_25}), $\frac{k}{\lambda}=0.15$~(\ref{boost-scale-fixed_0_15}), $\frac{k}{\lambda}=0.05$~(\ref{boost-scale-fixed_0_05}). The upper plot in each pair shows the zoomed version, focused on the region where both original and improved methods were able to converge.}\label{fig:boost-scale-fixed}
\end{figure}

\paragraph{Dependency on $\lambda$} We fix the threshold to the conservative value  $\frac{k}{\lambda} = 0.25$. We then plot the dependency of the run time on the non-centrality $\lambda$ for a set of d.o.f $k$: $2$, $15$, $50$, see Fig.~\ref{fig:boost-k-fixed}. For benchmarking we use the~\cite{gbenchmark} library. The benchmarking script itself became a part of the~\cite{gbench-script} library. Using this script we measure the run time 100 times and use the mean as a central value. The error bar is computed as a standard deviation. We add noise with standard deviation $\sigma = 10^{-6}$ to parameters $k$ and $\lambda$ to avoid caching effects. The vertical line on the plots shows the threshold where the original initial guess for smaller $\lambda$ is switched to the corrected value at bigger $\lambda$. Therefore, we expect that both lines coincide below the threshold and the improved solution would lie lower above the threshold. One can notice missing values on the curve representing the original initial guess. The reason for this is the numerical instability of the algorithm in Boost, that has been resolved after we corrected the initial guess.
\paragraph{Dependency on d.o.f ($k$)} In the set of plots in Fig.~\ref{fig:boost-scale-fixed}, we fix the asymptotic scale to values $\frac{k}{\lambda} = 0.25, 0.15, 0.05$ and investigate the dependency of the run time on the number of d.o.f. Since the threshold is fixed, the difference in the run time is caused by the actual position where the original initial guess starts to undershoot, the non-conservative threshold. The farther the fixed threshold is from the non-conservative threshold, the more significant the effect of undershooting at the test point will be. Therefore, we expect the difference in the run time to grow with number of d.o.f, as follows from Eq.~\ref{eq:boost-alpha-thresh}. For each value of the asymptotic scale, in addition to the full plot, we also show a zoomed version that shows the region where both original and improved methods were able to converge~\ref{fig:boost-scale-fixed}.

\section{Conclusion}
In this paper we present an approximate expression for the mode of the non-central $\chi^2$ distribution: $x_{mode} \approx k + \lambda - 3$, where $k$ is the number of degrees of freedom and $\lambda$ is the non-centrality parameter. The approximation is based on an asymptotic expansion and is valid in the region where the scale parameter $\frac{k}{\lambda} << 1$ and where the mode exists $k > 2$. The approximate formula can be used as the initial guess for iterative procedures searching for a precise solution. Run time performance and the domain of applicability of the Boost implementation of the mode search was improved using the presented approximate expression. The improvement became a part of the~\cite{gbench-script}.

\section{Acknowledgements}

We would like to thank Mykola~Semenyakin for the numerous fruitful and motivating discussions. We also would like to acknowledge the support of the Boost community that allowed the contribution to become a part of the Boost Library. This research was supported by the European Unions Framework Programme for Research and Innovation Horizon 2020 (2014-2021) under the Marie Sklodowska-Curie Grant Agreement No.765710.

\bibliographystyle{elsarticle-harv}
\bibliography{bibliography.bib}

\end{document}